\documentclass{article}

\usepackage[english]{babel}

\usepackage[letterpaper,top=2cm,bottom=2cm,left=3cm,right=3cm,marginparwidth=1.75cm]{geometry}

\usepackage{amsmath,amsthm}
\usepackage{graphicx}
\usepackage{paralist}
\usepackage{todonotes}

\usepackage{stex}
\usepackage{url}
\usepackage{booktabs}

\usepackage{enumitem}
\makeatletter
\def\namedlabel#1#2{\begingroup
    #2%
    \def\@currentlabel{#2}%
    \phantomsection\label{#1}\endgroup
}
\makeatother
\usepackage[colorlinks=true, allcolors=blue]{hyperref}

\title{How to Fight Fraudulent Publishing in the Mathematical Sciences: Joint Recommendations of the IMU and the ICIAM}

\author{Ilka Agricola, Lynn Heller, Wil Schilders, Moritz Schubotz, Peter Taylor, Luis Vega}

\date{September 2025}
\begin{document}
\maketitle

\textbf{ Prelude}.
\emph{In November 2023, Clarivate announced that it had excluded the entire field of
mathematics from the latest edition of its influential list of `highly cited researchers'. This prompted the IMU and the ICIAM to conduct a more thorough investigation into the
problem of fraudulent publishing in the mathematical sciences (see \cite{AHSSTV}). 
Understanding the problem is one thing; finding a way out and regaining control is another.
With the recommendations given below, we would like to start the discussion on how, as a global community, we can achieve this. We are all concerned.
 It affects the
very core of the science we love so much. I.A.}

\section{Introduction}
Predatory journals and citation cartels are reactions to the effort to exactly quantify and rank the quality of research through scientific `performance indicators', in the form of bibliometric measures. While evaluations are needed for the effective allocation of resources, a complete ranking of research quality cannot exist. Bibliometrics, on the contrary, suggest that everything can be ordered by simple numbers computed in some way or another from citations. In a field like mathematics, where the number of research papers and general citations are quite low, these citation numbers are prone to severe manipulation. The pressure to publish and the competitive job market provide an incentive for even serious scientists to artificially improve their own numbers. Predatory journals and citation cartels are at the extreme end of this enterprise, turning the desire to quantify into a money-making scheme. A fringe phenomenon in the past, this problem has by now reached even well-established research institutions and can no longer be ignored. 
Such manipulation results in a massive waste of limited resources being invested in low-quality research, while higher-quality work may be discontinued due to lack of funding. The opportunity cost of misaligned incentives is immense.

For a detailed description of fraudulent publishing in mathematical sciences, an explanation of technical terms,  and an extensive list of references, we refer to our article \cite{AHSSTV}. This note builds on that work to make explicit recommendations for policy makers and institutions on how to encourage
good scientific practice and simultaneously detect and combat fraudulent behavior. 
For individuals, we provide a set of possible actions to safeguard our research communities from being corrupted by fraud and pure business interests. 

These recommendations were formulated by the authors in close collaboration with the IMU Committee on Publishing (chaired by Ilka Agricola) and have been endorsed by the Executive Committee of the IMU and the Board of the ICIAM in May/June 2025. 


\section{Recommendations and Practical Advice}
The recent explosion in the number of publications in predatory journals demands action from all parties involved. The following list of recommendations for policy makers, institutions, and individuals maps out how we can all contribute to regain control over the situation. 

\subsection{For Policy Makers}
\subsubsection{Why should you care?}

\begin{itemize}

\item
Fraudulent publishing 
undermines trust in science and scientific results and therefore fuels anti-science movements. 

\item The use of falsified results can be dangerous and leads to a  waste of research effort and funding.


\item Fraudulent publication practice makes it hard for politicians and journalists to  distinguish serious science from junk science. 


\end{itemize}
\subsubsection{What can be done?}
\begin{itemize}
\item Assign resources based on expert-led assessments rather than bibliometric data. In particular, avoid relying on commercial journal rankings such as SJR and JCR.

\item Endorse good journals and discourage publishing in predatory journals. 

\item Define good publishing practices and encourage people to follow them.

\item Discourage the use of university rankings. 

\end{itemize}
\subsection{For Institutions}
\subsubsection{Why should you care?}

\begin{itemize}
\item 
Any instance of scientific fraud can cause significant long-term damage to an institution's reputation and result in a  failure to attract the best researchers and students and to raise grants.

\item Research evaluation based substantially on bibliometrics incentivizes manipulation, which can result in low quality scientists being hired or promoted. 



\item Fraud ruins scientific collaboration and morale and can mislead the younger generation 
into adopting bad publication habits.

\end{itemize}

\subsubsection{What can be done?}
\begin{itemize}
\item Discourage the use of bibliometrics in hiring and promotion committees.
\item Evaluate faculty by their best papers and `activeness'  without pressure to publish too frequently. 
\item  Educate researchers about predatory journals and discourage publication in these outlets.
\item Educate researchers about the low correlation between the quality of research and bibliometrics such as journal impact factors and citations.
\item Do not grant PhDs or other higher degrees based upon a requirement to publish a certain number of papers.
\item Choose carefully the Article Processing Charges (APCs) that you pay for. 
\item Implement best practice recommendations about the use of affiliations and professional email addresses.

\end{itemize}
\subsection{For Individuals}
\subsubsection{Why should you care?}
\begin{itemize}
    \item Your work  cannot compete with citation cartels if it is evaluated solely based on bibliometrics.
    \item Your  scientific integrity is at risk if you accidentally publish in a predatory journal.
    \item 
    If many in the community inflate their publication lists, the  pressure to increase output quantity will become even stronger.

\item Mathematics has many famous open  conjectures. 
Predatory journals can give people the opportunity to publish `proofs' of these without credible peer review. The status of these results can become unclear.

\end{itemize}
\subsubsection{What can be done?}
\begin{itemize}
    \item Read the actual publications 
    instead of relying on  bibliometrics, and say so when writing evaluations.  
    \item Avoid publishing in predatory journals
    (see Section \ref{predatory-journals}). 
    If you have done so in the past without knowing, 
    add a comment in your CV.
    \item Educate young researchers and colleagues about predatory publishing, the low correlation between bibliometrics and research quality, in particular in mathematics, and data literacy.
\item Help identifying good journals in your field; being listed by
zbMATH Open or MathSciNet is one good indicator of quality, but they will not cover some interdisciplinary outlets.

    \item Cite articles that are relevant for your work at hand, not more, not less.
    \item 
    If a respected colleague is associated with a predatory journal (as an author, reviewer, or guest editor), inform them of the journal's bad reputation and recommend that they stop cooperating with it.
    
    \item Check the quality of journals before joining their editorial board.
    \item Choose the journals or special issues in which you publish and the journals for which you write reviews wisely. 
    \item Be a responsible editor\,/\,editor-in-chief. Make any attempts at influencing you transparent. Resign if the situation becomes too bad, and make it public why you did it.
    \item What should you do if you find suspected misconduct or other irregularities in scientific papers? Make it public (but without putting yourself at risk).
   \item Be informed about screening tools for detecting plagiarism and fraud, and discuss their possible implementation if the need arises.
   \item Be critical! Verify author identities, affiliations, and email addresses whenever something seems suspicious. We recommend for this purpose ORCID (Open Researcher and Contributor ID), which is a unique and persistent identifier designed to solve the problem of name ambiguity. 
   

\end{itemize}

\subsubsection{Beware! }
\begin{itemize}
 \item The purpose of scientific phishing emails is not to trick the recipient into revealing sensitive information, but to draft them into the parallel universe of fake science at a low level. 

 \item Fraudulent journals often misuse the names of well-known researchers without their permission to create a false impression of credibility and legitimacy.

 \item Fraudulent journals often have bad publishing practices:
 not publishing accepted articles, 
  taking articles or journal websites offline without notice, publishing submitted articles before authors have signed a publishing agreement, and not reacting to emails.

\item `Hijackers' create web pages of existing, well-established journals to lure people into submitting to their fake journal; predatory journals often choose names that
are similar to those of existing prestigious journals.

\end{itemize}

\subsubsection{Strategic Recommendations: Alternatives needed} 
\begin{itemize}
\item We may need to create ``low entry" options of publication for papers that might not be of the highest level, but nevertheless are original and interesting.
\item Whenever you have the opportunity, advocate and fund \emph{serious science journalism}. It happens all too often that journalists stumble accross a non-reliable source and make a ``scoop" out of it.
\item
Governments should support efforts for greater transparency over fraudulent publications and
systematize large-scale monitoring of misconduct.
\item
A culture change in academia is needed that scientific misconduct is not acceptable at any level.

\end{itemize}

\section{Internet Resources}
Listed below are some Internet resources that provide valuable information on the topic.  
\subsection{Flagging Problematic Behavior}

\begin{compactdesc}
\item{\url{https://retractionwatch.com}} and its database, \url{http://www.retractiondatabase.org}/ 
 
  A free database of retracted scientific publications and an interesting news feature, now a part of \url{https://www.crossref.org}. If a journal or author accumulates retractions, be suspicious. 
\item{\url{https://pubpeer.com/}} Starting as a  “journal club” where members could discuss \emph{any} article post-publication, it has developed in recent years into a bulletin board for discussion of scientific misconduct, and serves sometimes as a starting point for further investigations by academic institutions or grant agencies. 
Although it is not well known in the mathematical community, it contains a certain amount of comments on mathematical papers, mostly from paper mills and predatory journals. 

\item \url{https://scienceintegritydigest.com/2019/07/16/how-to-report-misconduct-to-a-journal/}
This is a science integrity blog run by Elisabeth Bik. She won the Einstein Foundation Individual Award 2024 in recognition of her work, see
\url{https://award.einsteinfoundation.de/award-winners-finalists/recipients-2024/elisabeth-bik}

     \item{\url{https://clear-skies.co.uk}} Clear skies papermill alarm is a paid service that publishers can integrate into their peer review tools for detecting 
 papers that are close to ``past papermill-products", but there is also a 
 \href{https://clearskiesadam.medium.com/how-to-use-the-papermill-alarm-api-719b8b3b8253}{basic free API} for those with IT skills among us. It cooperates with the next tool.
\item{\url{https://www.irit.fr/~Guillaume.Cabanac/problematic-paper-screener}} The Problematic Paper Screener, a tool for flagging articles using tortured phrases.

\end{compactdesc}

\subsection{Checklists for identifying Predatory Journals}\label{predatory-journals}
\begin{compactdesc}
    \item{\url{https://thinkchecksubmit.org}}  A  non-profit cross-sector initiative aiming to
educate researchers, promote integrity, and build
trust in credible research and publications.
It provides guidelines, information, and free resources for recognizing predatory journals. 

\item Some universities or public institutions offer guidelines on how to identify predatory journals. We list here only a few examples.
\item{\url{https://researchguides.library.wisc.edu/c.php?g=1154500&p=8431826}} ---University of Wisconsin - Madison, US.
\item{\url{https://guides.libraries.indiana.edu/predatory}} ---
Indiana University, Bloomington, US.
\item{\url{https://www.centre-mersenne.org/editeurs-predateurs/}}
Centre Mersenne, France.

\end{compactdesc}

\subsection{Attempts for compiling lists of Predatory Journals}
It is virtually impossible to compile a complete list of predatory\,/\,junk journals.  For example, this `tag' may  not be permanent: there are examples of serious journals that became predatory and vice versa.

\begin{compactdesc}
\item{\url{https://beallslist.net/}} \ This was the first attempt to produce a list of potential predatory OA publishers, by the librarian Jeffrey Beall. As such, it had a tremendous influence on the perception of and discussions around predatory publishing. Several publishers sued Beall and his university because their journals were listed. The list has not been updated since 2017.

\item{\url{https://earlywarning.fenqubiao.com/#/en/early-warning-article-2024}} \ 
Since 2020, China has curated an `Early Warning Journal List' of predatory journals that researchers typically do not use anymore in subsequent years. It is not very long, and only ten journals with mathematical content have been listed so far (all of which can indeed safely be considered predatory). 
\item{\url{https://dgrsdt.dz/fr/revues_predateur}} Since 2018,
Algeria has also been compiling its own list of predatory journals and publishers. A quick look at the mathematical journals reveals that it is not complete, but those listed can again safely be considered predatory.
\item{\url{https://cabells.com/}} Cabell's List of predatory journals (paywalled). We cannot say anything about the quality of this list. 
\end{compactdesc}

\subsection{International Associations Promoting Scientific Integrity}
\begin{compactdesc}
    
\item{\url{https:// sfdora.org}}
The San Francisco Declaration on Research Assessment (DORA); its goal is  to develop new policies and practices for responsible research assessment. By now, their webpage also includes ``Reformscape", a searchable collection of criteria and standards for hiring, review, promotion, and tenure from academic institutions. 
All major science organizations have signed the declaration including the IMU, the ICIAM and the AMS.

\item{\url{https://publicationethics.org}  } The Committee On Publication Ethics (COPE) is an international organization with the mission to empower editors, journal staff, publishers, universities, organizations, and individuals by providing expert guidance, education and training, and offering a collaborative community. We recommend the Section
entitled \emph{Guidance}: ``Identify potential fraudulent activity and use the flowcharts when handling suspected manipulation. Suggested actions are recommended for each type of suspicious activity."

\end{compactdesc}

\subsection{Other resources}

\begin{compactdesc}

    \item{\url{https://doaj.org/}}\ This is a directory of Open Access journals, not more, not less. It awards a  ``DOAJ Seal" ``to journals that demonstrate best practice in open access publishing."
    To be honest, this is not a very hard quality criterion; keep away of OA journals that do not have the seal but do not take the seal as sole proof of quality. For example the MDPI journal \emph{Mathematics} has the seal. 
    Unfortunately, the list does not state any time period for which a journal has received the seal.

\item{\url{https://coara.eu}} \
In the future, the COARA initiative\footnote{\url{https://coara.eu/app/uploads/2022/09/2022_07_19_rra_agreement_final.pdf}} is planning to establish other metrics, which will increase efforts to measure different kinds of scientific contributions. It is not yet clear whether this will improve or worsen the situation.  
\item{\url{https://papermills.tilda.ws/advent}}\
This is a special `Advent calendar' designed by science sleuth Anna Albakina. Readers are asked to identify problems and inconsistencies in twenty-four cases involving academic papers submitted to scientific journals. 

 \end{compactdesc}

\end{document}